\documentclass[secthm,seceqn,amsthm,ussrhead,eqno]{amsart}
\usepackage{amsmath,latexsym}
\usepackage[psamsfonts]{amssymb}
\usepackage{times}
\usepackage[mathcal]{euscript}
\numberwithin{equation}{section}

\newcommand{\bbT}{\mathbb T}

\renewcommand{\epsilon}{\varepsilon}

\newcommand{\be}{\begin{equation}}
\newcommand{\ee}{\end{equation}}
\newcommand{\no}{\nonumber}

\newcommand{\C}{\mathbb{C}}

\newcommand{\N}{\mathbb{N}}

\newcommand{\R}{\mathbb{R}}
\renewcommand{\S}{\mathbb{S}}
\newcommand{\T}{\mathbb{T}}

\newcommand{\Z}{\mathbb{Z}}

\newcommand{\cB}{{\mathcal B}}

\renewcommand{\det}{\mathop{\mathrm{det}}}

{\bf}{\it}
\newtheorem{theorem}{Theorem}[section]
\newtheorem{lemma}[theorem]{Lemma}

\newtheorem{corollary}[theorem]{Corollary}
\newtheorem{assumption}[theorem]{Assumption}

\newtheorem{definition}[theorem]{Definition}

\newtheorem{remark}[theorem]{Remark}

\date{\today}

\begin{document}
\title[The threshold effects for a
family of Friedrichs models...] {The threshold effects for a family
of Friedrichs models under rank one perturbations}

\author{Sergio  Albeverio$^{1,2,3}$, Saidakhmat  N.~Lakaev$^{4,5}$,
  Zahriddin I.~Muminov $^{5}$}

\address{$^1$ Institut f\"{u}r Angewandte Mathematik,
Universit\"{a}t Bonn, Wegelerstr. 6, D-53115 Bonn\ (Germany)}

\address{
$^2$ \ SFB 611, \ Bonn, \ BiBoS, Bielefeld - Bonn;}
\address{
$^3$ \ CERFIM, Locarno and Acc.ARch,USI (Switzerland) E-mail:
albeverio@uni.bonn.de}

\address{
{$^4$ Samarkand State University, University Boulevard 15, 703004
Samarkand (Uzbekistan)} \ {E-mail: lakaev@yahoo.com }}

\address{$^5$ Samarkand division of Academy of sciences of
Uzbekistan (Uzbekistan) E-mail:~zimuminov@mail.ru }

\begin{abstract}
A family of  Friedrichs models under rank one perturbations
$h_{\mu}(p),$  $p \in (-\pi,\pi]^3$, $\mu>0,$ associated to a system
of two particles on the three dimensional lattice $\Z^3$ is
considered. We prove the existence of a unique eigenvalue below the
bottom of the essential spectrum of $h_\mu(p)$ for all nontrivial
values of $p$ under the assumption that $h_\mu(0)$ has either a
threshold energy resonance (virtual level) or a threshold
eigenvalue. The threshold energy expansion for the Fredholm
determinant associated to a family of Friedrichs models is also
obtained.
\end{abstract}
 \maketitle



\section{Introduction}
In the present paper we consider a family of Friedrichs models under
rank one perturbations associated to a system of two particles on
the lattice $\Z^3$ interacting via pair non-local potentials.

The main goal of the paper is to give a thorough mathematical
treatment of the spectral properties of a family of Friedrichs
models in dimension three with emphasis on threshold energy
expansions for the Fredholm determinant associated to the family
(see, e.g. \cite{AGH,AGHH,ALMM, ALzMahp,FIC,Lfa93,Sob,Tam94,Yaf2}
for relevant discussions and \cite{K-M,M-S,Zh} for the general study
of the low-lying excitation spectrum for quantum systems on
lattices).

These kind of models have been discussed in quantum mechanics
\cite{Fad,Fri}, solid physics \cite{RSIII, Mat, Mog,G-Sh} and in
lattice field theory \cite{MalMin,LMfa04,LMfa05}.

Threshold energy resonances (virtual levels) for the two-particle
Schr\"odinger operators have been studied in
\cite{AGH,ALzMahp,ALMM,K-S,Lfa93,Yaf2}. Threshold energy expansions
for the resolvent of two-particle Schr\"odinger operators have been
studied in \cite{ALzMahp,J-K,K-S,Lfa93,Ltmf92,Sob,Tam94,Yaf2} and
have been applied to the proof of the existence of Efimov's effect
in \cite{ALzMahp,Lfa93,Sob,Tam94,Yaf1}.

Similarly to the lattice  Schr\"odinger operators and in contrast to
the continuous Schr\"o\-din\-ger operators the family of Friedrichs
models $h_{\mu}(p),$  $p \in (-\pi,\pi]^3$, $\mu>0$ depends
parametrically on the internal binding $p $, the quasi-momentum,
which ranges over a cell of the dual lattice and hence it has
spectral properties analogous to those of lattice Schr\"odinger
operators.

Let us recall that the spectrum and resonances of the original
Friedrichs model and its generalizations have been studied  and the
finiteness of the eigenvalues lying below the bottom of the
essential spectrum has been proven in \cite{Fri,Fad,Ltsp86,Yaf3}.

In \cite{LMfa04,LMfa05} a peculiar  family of Friedrichs models was
considered and the appearance of eigenvalues for values of the total
quasi-momentum $p \in (-\pi,\pi]^d,d=1,2$ of the system lying in a
neighbourhood of some particular values of the parameter $p$ has
been proven.

 For a wide class of the two-particle Schr\"odinger operators only the
existence of eigenvalues of $h_\mu(p), p \in (-\pi,\pi]^3,$ for all
nonzero values of the quasi-momentum $0\neq p\in \T^3$ (under the
assumption that  $h(0)$ has either a zero energy resonance or a zero
eigenvalue) has been proven in \cite{ALMM}.

In the present paper two main results.

First of them gives the existence of a unique eigenvalue $e_\mu(p)$
of $h_\mu(p), p \in (-\pi,\pi]^3,$ for all nonzero values of the
quasi-momentum $0\neq p\in \T^3$ (provided that $h_\mu(0)$ has
either a threshold  energy resonance or a threshold eigenvalue) and
lower and upper bounds on it. The monotonous dependence of the
eigenvalue $e_\mu(p)$ on $\mu$ (Theorem \ref{pos.eig}).

The second one presents an expansion for the Fredholm determinant
resp. the Birman-Schwinger operator in powers of the quasi-momentum
$p \in (-\pi,\pi]^3$ in a small $\delta$-neigh\-bor\-hood of the
origin and proving that the Fredholm determinant resp. the
Birman-Schwinger operator has a differentiable continuation to the
bottom of the essential spectrum of $h_{\mu}(p)$ as a function of
$w=(m-z)^{1/2}\geq 0$ for $z\leq m,$ where $z\in \R^1$ is the
spectral parameter (Theorem \ref{main}).

We notice that if the functions $u$ resp. $\varphi$ is analytic on
$(\T^3)^2$ resp. $\T^3,$ then one can obtain a precise expansion for
the Fredholm determinant and the Birman-Schwinger operator (see
\cite{K-S,Ltmf92}).

 The structure of the
paper is as follows. In Section 2 we state the problem and present
the main results. Proofs are presented in Section 4 and are based on
a series of lemmas in Section 3. In Appendix for an important
subclass of the family of Friedrichs models we shall show that all
assumptions in Section 2 are fulfilled.

Throughout the present paper we adopt the following conventions:
$\T^3$ denotes  the three-dimensional torus, the cube $(-\pi,\pi]^3$
with appropriately  identified sides. For each $\delta>0$ the
notation $U_{\delta}(0) =\{p\in {\T}^3:|p|<\delta \}$ stands for a
sufficiently small $\delta$-neighborhood of the origin. Denote by
$L_2(\T^3)$ the Hilbert space of square-integrable functions on
$\T^3.$

Let  $\cB (\theta,U_{\delta}(0))$ with $1/2<\theta<1$,  be the
Banach spaces of H\"older continuous functions on
$\overline{U_{\delta}(0)}$ with exponent $\theta$ obtained by the
closure of the space of smooth functions $f$ on $U_{\delta}(0)$ with
respect to the norm
\begin{equation*}
\|f\|_{\theta}=\sup_{t, \ell\in U_{\delta}(0) \atop t\neq \ell
}\bigg [|f(t)|+|t-\ell|^{-\theta}|f(t)-f(\ell)|\bigg ].
\end{equation*}
The set of functions $f: {\bbT}^3\to \R$ having  continuous partial
derivatives up to order $\leq n$ will be denoted by
$C^{(n)}({\bbT}^3).$ In particular $C^{(0)}({\bbT}^3)=C({\bbT}^3).$

\section{The  model operator $h_\mu(p)$, main assumptions and
statement of the results}

Let $u$ be a real-valued
 essentially bounded function on $({\T}^3)^2$ and $\varphi$ be
a real-valued function in $L_2(\T^3).$ Let $\mu$ be a positive real
number.

We introduce the following family of bounded self-adjoint operators
(the Friedrichs model) $h_\mu(p),\,  p\in {\T}^3,$ acting in
$L_2(\T^3)$ by
\begin{equation*}\label{h_alpha}
 h_\mu(p)=h_{0}(p)-\mu v,
 \end{equation*}
 where
\begin{align*}
 &(h_0(p)f)(q)=u(q,p)f(q),\,\,
 f\in L_2(\T^3),\\
\end{align*} and $v$ is non-local interaction operator
\begin{align*}\label{poten}
 &(vf)(q)=\varphi(q)\int_{\T^3}
\varphi(t)f(t)dt,\,\,
 f\in L_2(\T^3).
\end{align*}

\begin{remark}
In the case where the function $u$ is of the form
\begin{equation}\label{of the form}
u(p,q)=\varepsilon (p)+ \varepsilon (p-q)+ \varepsilon (q),
\end{equation}
 the operator $h_0(p)$ is associated
to a system of two particles (bosons) moving on the
three-dimensional lattice $\Z^3$ and is called the {\it\bf  free
Hamiltonian}, where $\varepsilon (\cdot)$ is  the {\it dispersion
relations of normal modes}  associated with  the free particle in
question.
\end{remark}

Throughout this paper we assume  the following additional
hypotheses.

\begin{assumption}\label{hypoth} $(i)$ The function $u$ is even on
$({\T }^3)^2$ with respect to $(p,q),$ and has a unique
non-degenerate
 minimum at the point $(0,0)\in ({\T}^3)^2$ and all third order
partial derivatives  of  $u$ are continuous on $(\T^3)^2$ and their
restrictions in $(U_{\delta}(0))^2$ belong to
$\cB (\theta,(U_{\delta}(0))^2).$\\
$(ii)$ For some positive definite matrix $U$ and real numbers
$l,l_1, l_2 \,(l_1,l_2>0,l\not=0)$ the following equalities hold
$$
\left( \frac{\partial^2 u(0,0)}{\partial p^{(i)} \partial p^{(j)}}
\right)_{i,j=1}^3= l_1 U,\,\, \left( \frac{\partial^2
u(0,0)}{\partial p^{(i)} \partial q^{(j)}} \right)_{i,j=1}^3= l
U,\,\, \left( \frac{\partial^2 u(0,0)}{\partial q^{(i)} \partial
q^{(j)}} \right)_{i,j=1}^3= l_2 U.
$$
\end{assumption}

\begin{remark}
The function $u$ is even and has a unique non-degenerate minimum on
$\T^3$ and hence without loss of generality we assume that the
function $u$ has a unique minimum at the point $(0,0)\in (\T^3)^2.$
\end{remark}
\begin{remark}
It is easy to see that Assumption \ref{hypoth} implies the
inequality \, $l_1l_2 > l^2.$
\end{remark}

\begin{assumption}\label{hyp.varphi} The continuous function
 $\varphi$  is either even or odd on
$\T^3$ and all second order partial derivatives  of   $\varphi$ are
continuous on $\T^3.$
\end{assumption}

Let ${\C}$ be the field of complex numbers and  set
$$ u_p(q)=u(p,q),\quad
\quad m=\min_{p,q\in \T^3} u(p,q),
  $$\begin{equation*}
u_{\min}(p)=\min_{q\in {\T}^3}u_p(q),\quad u_{\min}(p)=\max_{q\in
{\T}^3} u_p(q)
\end{equation*}
 and
\begin{equation}\label{Lamb}
\Lambda(p,z)=\int\limits_{{\T}^3}
\frac{\varphi^2(t)dt}{u_p(t)-z},\quad p\in\T^3,\,\, z \in
\C\setminus [u_{\min}(p),u_{\max}(p)].
\end{equation}

\begin{remark}\label{belongs}
By part $(i)$ of Assumption \ref{hypoth} all third order
partial derivatives of the function  $\Lambda(\cdot,z),$ $z< m,$
 belong to
$C^{(2)}({\T}^3).$
\end{remark}

 The function
$\Lambda(p,\cdot),p\in \T^3,$ is increasing in
$(-\infty,u_{\min}(p))$ and hence the following finite or infinite
positive limit exists
\begin{equation}\label{limm}
    \lim_{z\to u_{\min}(p)-0}\Lambda(p,z)=\Lambda(p,u_{\min}(p)).
\end{equation}

\begin{remark}\label{limit}
 Since for any $p\in U_\delta(0),$ $\delta>0$-sufficiently small, the
function $u_p(\cdot)$ has a unique non-degenerate minimum in $\T^3$
(see part $(i)$ of Lemma \ref{minimum}) Lebesgue's dominated
convergence theorem  yields the equality
\begin{equation*}
\Lambda(p,u_{\min}(p))=\int\limits_{{\T}^3}
\frac{\varphi^2(t)dt}{u_p(t)-u_{\min}(p)},\quad p\in U_\delta(0).
\end{equation*}
\end{remark}

 The perturbation $v$ of the multiplication operator
$h_0(p)$ is of rank one and hence, in accordance with  Weyl's
theorem the essential
  spectrum of the operator
$h_{\mu}(p)$ fills the following interval on the real axis:
$$
\sigma_{ess}(h_{\mu}(p))=[u_{\min}(p),u_{\max}(p)].
$$

\begin{remark}
We remark that for some $p\in \T^3$  the essential spectrum of
$h_{\mu}(p)$ may degenerate to the set consisting of the unique
point $[u_{\min}(p),u_{\min}(p)].$ Because of this we can not state
that the essential spectrum of $h_{\mu}(p)$ is absolutely continuous
for any $p\in \T^3.$ This is the case, e.g. for a function $u$ of
the form \eqref{of the form}, where $p=(\pi,\pi,\pi)\in \T^3,$ and
\begin{align} \label{u of forms} &\varepsilon (q)=3-\cos q_1-\cos
q_2 -\cos q_3, \quad q=(q_1,q_2,q_3) \in {\T}^3.
\end{align}

\end{remark}

\begin{definition}\label{resonance0} Assume part $(i)$ of Assumption
\ref{hypoth} and $ \varphi \in {C^{(0)}(\T^3)}.$ The operator
$h_\mu(0)$ is said to have a threshold energy resonance if the
number  $1$ is an eigenvalue of the operator
$$
(\mathrm{G}\psi)(q)=\mu \varphi(q) \int\limits_{{\T}^3}
\frac{\varphi(t)\psi(t)dt} {u_0(t)-m},\quad \varphi\in
{C^{(0)}(\T^3)}
$$
and the associated eigenfunction $\psi $ (up to constant factor)
satisfies the condition $\psi(0)\neq 0.$
\end{definition}

\begin{remark}
Assume part $(i)$ of Assumption \ref{hypoth} and Assumption
\ref{hyp.varphi}.\\
(i) If $\varphi(0)\not=0$ and the operator $h_{\mu}(0)$ has a
threshold energy resonance, then the function
\begin{equation}\label{f0f1}
f(q)=\frac{\varphi(q)}{u_0(q)-m},
\end{equation}
obeys the equation $ h_{\mu}(0)f=mf$ and  $ f\in L_1(\T^3)\setminus
L_2(\T^3)$
(see Lemma \ref{resonance}).\\
(ii) If $\varphi(0)=0$ and the threshold $z=m$ is an eigenvalue of
the operator $h_{\mu}(0),$ then the function $f,$ defined by
\eqref{f0f1}, obeys the equation $ h_{\mu}(0)f=mf$ and $ f\in
L_2(\T^3)$ (see Lemma \ref{zeroeigen}).
\end{remark}

Set
\begin{equation*}\label{mu.alpha}
 \mu_0=\Lambda^{-1}(0,m).
\end{equation*}

\begin{remark} Notice that the conditions $\mu =\mu_0$ and
$\varphi(0)\neq 0$ (resp. $\mu =\mu_0$ and $\varphi(0)=0$) mean that
the operator $h_{\mu}(0)$ has threshold energy resonance (see Lemma
\ref{resonance}) (resp.
  a threshold eigenvalue of $h_{\mu}(0)$ (see Lemma \ref{zeroeigen})).
\end{remark}

\begin{remark}\label{pos.eig}
We note that the  bottom $z=m$ of the essential spectrum
$\sigma_{ess}(h_{\mu_0}(0))$ of $h_{\mu_0}(0)$ is either a threshold
energy resonance or an eigenvalue for the operator $h_{\mu_0}(0).$
\end{remark}

In order to  study the spectral properties of $h_{\mu}(p)$ precisely
we assume the following
\begin{assumption}\label{Lambda}
Assume that: (i) the function $\Lambda(\cdot,u_{min}(\cdot))$ has a
unique minimum at the origin, i.e. for all $0\neq p\in \T^3$ the
following inequality holds
\begin{equation*}\label{requested}
\Lambda(p,u_{\min}(p))-\Lambda(0,u_{\min}(0))>0.
\end{equation*}
$(ii)$ the function $\Lambda(\cdot,m)$ has a unique maximum at the
origin such that for some $c>0$ the following inequality holds
\begin{equation*}\label{}
\Lambda(0,m)-\Lambda(p,m)>c|p|^2,\quad 0\neq p\in U_\delta(0).
\end{equation*}

\end{assumption}
\begin{remark}
If for all $0\neq p\in \T^3$ and a.e. $ q \in \T^3$ the inequality
$$u_p(q)-u_{\min}(p)<u_0(q)-u_{\min}(0)$$
holds, then part (i) of Assumption \ref{Lambda} is obviously
fulfilled. In  Appendix we shall show that for the functions of the
form \eqref{of the form}   Assumption \ref{Lambda} is fulfilled.
\end{remark}
The following theorem presents a result characteristic for the
two-particle Hamiltonians on lattices (see \cite{ALMM}).

\begin{theorem}\label{pos.eig}
Assume Assumptions \ref{hypoth}, \ref{hyp.varphi} and \ref{Lambda}.
Then for all $p\in \T^3\setminus \{0\}$ the operator $h_{\mu_0}(p)$
has a unique eigenvalue $e_{\mu_0}(p).$ One has
$$m<e_{\mu_0}(p)<u_{\min}(p),\,0\neq p\in
\T^3.$$ (ii)  For any $\mu>\mu_0$ the operator $h_\mu(p),\,p\in
\T^3$ has a unique eigenvalue $e_\mu(p).$ One has
$$e_\mu(p)<e_{\mu_0}(p)<u_{\min}(p),\,0\neq p\in
\T^3$$ and
$$e_\mu(0)<m.$$
\end{theorem}

 For any $p \in \T^3$ we
define an analytic  function $\Delta_{\mu}(p,\cdot)$ (the Fredholm
determinant
 associated to the operator $h_\mu(p)$)  in
 ${\C} { \setminus } [u_{\min}(p),u_{\max}(p)]$ by
\begin{equation*}\label{det}
\Delta_{\mu}(p,z)=1-\mu
\int\limits_{{\T}^3}(u(p,t)-z)^{-1}\varphi^2(t)dt.
\end{equation*}

Now we formulate a result (threshold energy expansion for the
Fredholm determinant) of the paper, which is important in the
spectral analysis for a model operator associated to a system of
three-particles on the lattice $\Z^3$ \cite{ALZMarX06}.
\begin{theorem}\label{main} Assume Assumptions \ref{hypoth} and \ref{hyp.varphi}.

 For any  $z< u_{\min}(p)$ the function $\Delta_\mu(\cdot,z)$
 is of class
$C^{(2)}(\T^3)$  and the following decomposition
\begin{align*}
& \Delta_\mu(p,z)=\Delta_\mu(0,m)+\\
 &\frac{2\sqrt{2}\pi^2 \mu
\varphi^2(0)}
 {l_1^{\frac{3}{2}} \det(U)^{\frac{1}{2}}}\sqrt{u_{min}(p)-z}+
\Delta_\mu^{res}(u_{min}(p)-z)+\Delta_\mu^{res}(p,z),\,\,z\leq
u_{min}(p),\,\,p\in U_\delta(0),
\end{align*}
holds, where
$\Delta_\mu^{res}(u_{min}(p)-z)=O((u_{min}(p)-z)^{(1+\theta)/2})$ as
$z \to u_{min}(p)$ $(z<u_{\min}(p))$ and
 $\Delta_\mu^{res}(p,z)=O(p^2)$ as $p\to 0$ uniformly in
$z\leq u_{\min}(p).$
\end{theorem}
\begin{remark}
An analogue result has been proven in \cite{ALzMahp} in the case
 where $\varphi(\cdot)\equiv const$ and the function $u(\cdot,\cdot)$
is of the form \eqref{of the form}, \eqref{u of forms}.
\end{remark}

\begin{corollary}\label{razl.lemma.natijasi.}
Assume Assumptions \ref{hypoth} and \ref{hyp.varphi}.

 (i) Let the operator $h_{\mu_0}(0)$ have a thrshold energy
resonance. Then for all $p\in U_{\delta}(0) $ and $z \leq m$ the
following decomposition holds
\begin{align*}
\Delta_{\mu_0}(p,z)=\frac{4\sqrt{2}\pi^2 \mu_0 \varphi^2(0)}
 {l_1^{\frac{3}{2}} \det(U)^{\frac{1}{2}}}\sqrt{u_{min}(p)-z}+
\Delta_{\mu_0}^{res}(u_{min}(p)-z)+ \Delta_{\mu_0}^{res}(p,z).
\end{align*}

(ii) Let the threshold $z=m$ be an  eigenvalue of $h_{\mu_0}(0)$.
Then for any $p\in U_{\delta}(0) $ and $z \leq m$ the following
decomposition holds
\begin{align*} \Delta_{\mu_0}(p,z)= \Delta_{\mu_0}^{res}(u_{min}(p)-z)+
\Delta_{\mu_0}^{res}(p,z).
\end{align*}
\end{corollary}
\begin{remark}
We see that Corollary \ref{razl.lemma.natijasi.} gives a threshold
energy expansions for the Fredholm determinant, leading to different
behaviors for a threshold energy resonance resp. a threshold
eigenvalue.
\end{remark}
The following Corollary \ref{D.ineq} (resp. Corollary
\ref{fin.ineq}) plays a crucial role in the proof of the
infiniteness (resp. finiteness) of the number of eigenvalues lying
below the bottom of the essential spectrum for a model operator
associated to a system of three-particles on  three dimensional
lattice $\Z^3$ \cite{ALZMarX06}.
\begin{corollary}\label{D.ineq} Assume Assumptions \ref{hypoth} and \ref{hyp.varphi}.
Let the operator $h_{\mu_0}(0)$
 have a threshold energy resonance.
 Then for some
$c_1,c_2>0$  the following inequalities hold
\begin{equation}\label{c<(.,.)<c}
c_1 |p| \leq \Delta_{\mu_0}(p,m) \leq c_2 |p|,\quad p\in
U_\delta(0),
 \end{equation}
\begin{equation}\label{(.,.)>c}
\Delta_{\mu_0}(p,m) \geq c_1,  \quad p\in \T^3\setminus U_\delta(0).
 \end{equation}

\end{corollary}

\begin{corollary}\label{fin.ineq} Assume Assumptions \ref{hypoth},
\,\ref{hyp.varphi} and
\ref{Lambda}. Let $z=m$ be an  eigenvalue of $h_{\mu_0}(0)$. Then
for some $c>0$ the following inequality holds
\begin{align*}
\Delta_{\mu_0}(p,m)\geq c p^2,\quad p\in U_\delta(0).
\end{align*}

\end{corollary}

\section{Spectral properties of the operator $h_\mu(p)$}

In this section we study some spectral properties of the family
$h_\mu(p),\,  p\in {\T}^3,$ with emphasis on the threshold energy
resonance and a threshold eigenvalue.

Denote by $r_{0}(p,z)=(h_0( p)-zI)^{-1},$ where $I$ is identity
operator on $L_2(\T^3),$ the resolvent of the operator $h_0(p),$
that is, the multiplication operator by the function
$$(u_p(\cdot)-z)^{-1},\quad z\in \C\setminus [u_{\min}(p),u_{\max}(p)].$$

From the equality \eqref{limm} we have
$$
\Delta_{\mu}(p,m)=\lim_{z\to m-0} \Delta_{\mu}(p,z),\quad  p\in
U_\delta(0).
$$

\begin{lemma}\label{delta=0}
 Assume part $(i)$ of Assumption
\ref{hypoth} and $ \varphi \in {C^{(0)}(\T^3)}.$ For any $\mu>0$ and
$p\in \T^3$ the following statements are equivalent:\\
(i) the operator $h_\mu(p)$ has an eigenvalue $z \in {\C} \setminus
[u_{\min}(p),u_{\max}(p)]$ below the bottom of the essential
spectrum.\\
(ii) $\Delta_{\mu}(p,z)=0$, $z \in
{\C} \setminus [u_{\min}(p),u_{\max}(p)].$\\
 (iii) $\Delta_{\mu}(p,z')<0$ for some
$z'\leq u_{\min}(p).$
\end{lemma}
\begin{proof}
We prove\,\,$(i)\rightarrow(ii)\rightarrow(iii)\rightarrow(ii).$
From the positivity of $v$ it follows that the positive square root
of $v$ exists, we shall denote it by $v^{\frac{1}{2}}.$ For any
$\mu>0$ and $p\in \T^3$ the number $z \in {\C} \setminus
[u_{\min}(p),u_{\max}(p)]$ is an eigenvalue of $h_\mu(p)$ if and
only if $\lambda=1$ is an eigenvalue of the operator
\begin{equation*}
G_{\mu}(p,z)=\mu v^{\frac{1}{2}} r_0(p,z)v^{\frac{1}{2}}
\end{equation*}
(this follows from the Birman-Schwinger principle).

Since the operator $v^{\frac{1}{2}} $ is of the form $$
(v^{\frac{1}{2}} f)(q)=||\varphi||^{-1}\varphi(q)\int\limits_{\T^3}
\varphi(t)f(t)dt,\quad f \in L_2(\T^3) $$ the operator
$G_{\mu}(p,z)$ has the form
\begin{equation*}
(G_{\mu}(p,z)f)(q)=\frac{\mu {\Lambda}(p,z) }{||\varphi||^2}
\varphi(q)\int\limits_{\T^3} \varphi(t)f(t)dt,\quad f \in L_2(\T^3),
\end{equation*}
where $\Lambda(p,z)$ is defined by \eqref{Lamb}.

According to the Fredholm theorem  the number $\lambda=1$ is an
eigenvalue of the operator $G_{\mu}(p,z)$ if and only if
$\Delta_{\mu}(p,z)=0.$  $(i) \Leftrightarrow (ii)$ is proven.

$(ii) \Leftrightarrow (iii)$.
 Let $\Delta_{\mu}(p,z_0)=0$ for some $z_0 \in
\C\setminus  [u_{\min}(p),u_{\max}(p)]$. The operator $h_\mu(p)$ is
self-adjoint and
 hence by $(i)\Leftrightarrow (ii)$ we conclude that $z_0$ is a real.
For all $z>u_{\max}(p)$ we have $\Delta_{\mu}(p,z)>1$ and hence
$z_0\in (-\infty,u_{\min}(p)).$
 Since for any $p\in
{\T}^3$ the function $\Delta_{\mu}(p,\cdot)$ is decreasing in $z\in
(-\infty,u_{\min}(p))$ we have $\Delta_{\mu}(p,z')<
\Delta_{\mu}(p,z_0)=0$ for some $ z_0<z'\leq u_{\min}(p).$

$(iii) \Leftrightarrow (ii)$. Suppose that $\Delta_{\mu}(p,z')<0$
for some $z' \leq u_{\min}(p).$ For any $p\in {\T}^3$ we have
$\lim\limits_{z\to -\infty} \Delta_{\mu}(p,z)=1,$
$\Delta_{\mu}(p,\cdot)$ is continuous in $z\in
(-\infty,u_{\min}(p))$ and hence there exists $z_0\in (-\infty,z')$
such that $\Delta_{\mu}(p,z_0)=0.$ This completes the proof.
\end{proof}
The following lemmas  describe whether the bottom of the essential
spectrum of $h_{\mu_0}(0)$ is threshold energy resonance or a
threshold eigenvalue.

\begin{lemma}\label{resonance} Assume part $(i)$ of Assumption \ref{hypoth} and
 $ \varphi \in {\mathcal{B}(\theta,\T^3)} ,\,\,\frac{1}{2}<\theta\leq 1.$
The following statements are equivalent:\\
(i) the operator $h_{\mu}(0)$ has a threshold energy resonance and
\begin{equation}\label{resfunction}
\varphi(q)(u_0(q)-m)^{-1} \in L_1(\T^3)\setminus
L_2(\T^3).
\end{equation}
(ii) $\varphi(0)\neq 0$ and
$\Delta_{\mu}(0,m)=0.$\\
(iii) $\varphi(0)\neq 0$  and $\mu= \mu_0.$
\end{lemma}

  \begin{proof} We
prove\,\,$(i)\rightarrow(ii)\rightarrow(iii)\rightarrow(i).$ Let the
operator $h_{\mu}(0)$ have threshold energy resonance  for some
$\mu>0$. Then by
 Definition \ref{resonance0} the equation
\begin{equation}\label{res-def}
\psi(q)=(G \psi)(q),\quad \psi\in C(\T^3)
 \end{equation}
 has a simple solution $\psi(\cdot)$ in $C({\T^{3}}),$ such that $\psi(0)\neq
 0.$

This solution is equal to the  function $\varphi$ (up to a constant
factor) and hence $ \Delta_{\mu}(0,m)=0$ and so $\mu=\mu_0$.

 Let $\varphi(0)\neq 0$ and $\mu=\mu_0,$ hence  the equality $ \Delta_{\mu}(0,m)=0$
 holds.
Then  the function $\varphi$ is a solution of  equation
\eqref{res-def}, that is, the operator $h_{\mu}(0)$ has a threshold
energy resonance. Since $u_0(\cdot)$ has a unique non-degenerate
minimum at $p=0\in \T^3$ and $\varphi(0)\neq 0$ the inclusion
\eqref{resfunction} holds.
\end{proof}

\begin{lemma}\label{zeroeigen} Assume part $(i)$ of Assumption \ref{hypoth}
and $ \varphi \in {\mathcal{B}(\theta,\T^3)},\,\,
\frac{1}{2}<\theta\leq 1.$
The following statements are equivalent:\\
(i) the operator $h_{\mu}(0)$ has    a threshold eigenvalue.\\
(ii) $\varphi(0)=0$ and $\Delta_{\mu}(0,m)=0.$ \\
(iii) $\varphi(0)=0$ and $\mu= \mu_0.$
\end{lemma}
\begin{proof}
We prove\,\,$(i)\rightarrow(ii)\rightarrow(iii)\rightarrow(i).$
Suppose $f\in L_2(\T^3)$ is an eigenfunction of the operator
$h_{\mu}(0)$ associated with the  eigenvalue $m$. Then
\begin{equation}\label{h=0}
(u_0(q)-m)f(q)-\mu \varphi(q) \int\limits_{{\T}^3}
\varphi(t)f(t)dt=0  \quad\mbox {and}\quad \int\limits_{{\T}^3}
\varphi(t)f(t)dt\neq 0.
\end{equation}
Hence \eqref{h=0} yields $\varphi(0)=0.$ We find that $f,$ except
for an arbitrary factor, is given by
\begin{equation}\label{eigenunction}
f(q)=(u_0(q)-m)^{-1}\varphi(q).
\end{equation}
Thus \eqref{h=0} implies $\Delta_{\mu}(0,m)=0$ and $\mu=\mu_0.$

Let $\varphi(0)=0$ and  $\mu=\mu_0,$ then
 $\Delta_{\mu}(0,m)=0$ and the function
$f$, defined by \eqref{eigenunction}, obeys the equation
$h_{\mu}(0)f=m f.$ Since $u_0(\cdot)$ has a unique non-degenerate
minimum at $p=0\in \T^3$ and $\varphi(0)=0$ we have
 $ f\in L_2(\T^3).$
 \end{proof}

\begin{lemma}\label{minimum} Assume Assumption \ref{hypoth}  be fulfilled.
Then there exists a $ {
\delta } $-neighborhood $U_ {\delta } (0)\subset \T^3$ of the point
$p=0$ and a function $q_0 (\cdot)\in C^{(2)}(U_ {\delta
} (0))$   such that:\\
(i)   for any $p { \in } U_ { \delta }(0)$  the point
  $q_0(p) $ is a unique non-degenerate minimum  of $u_p(\cdot)$ and
\begin{equation}\label{q0 asym}
q_0(p)=-\frac{l_2}{l_1}p+O(|p|^{2+\theta})\,\,as\,\,p \to 0.
\end{equation}\\
(ii) the function $u_{min}(\cdot)= U(\cdot,q_0(\cdot))$ is even,  of
class $ C^{(3)}(U_ {\delta } (0))$  and has the asymptotics
 \begin{equation}\label{min.raz}
u_{min}(p)=m+\frac{l^2_1-l^2_2}{2l_1}(Up,p)+O(|p|^{3+\theta}) \quad
\mbox{as}\quad p \to 0.
\end{equation}
(iii) let for some $p\in \T^3$ the point $q_0(p)$  be a minimum  of
$u_p(\cdot)$ (if the minimum value of
$\hat{u}_{p}(\cdot)=\hat{u}(p,\cdot)$
 is attained in several points
$q_0(p)$ as nearest to $0\in T^3$), that is, $u_p(q_0(p))=\min_{q\in
\T^3} u_p(q)$. Then $q_0(-p)=-q_0(p).$
\end{lemma}
\begin{proof}
$(i)$ By the implicit function theorem there exist ${\delta }>0$ and
a function $q_0 (\cdot)\in C^{(2)}(U_ {\delta } (0))$  such that for
any $p { \in } U_ { \delta} ( 0 ) $ the point $ q_0(p)$ is the
unique non-degenerate minimum point of $u_p(\cdot)$ (see Lemma 3 in
\cite{Ltmf92}).

Since $u(\cdot,\cdot)$ is even with respect to $(p,q)\in (\T^3)^2$
we have
\begin{equation}\label{evenn}
u_{min}(-p)= \min_{q\in \T^3} u_{-p}(q) = \min_{q\in \T^3} u_p(-q)=
\min_{-q\in \T^3} u_p(q)= \min_{q\in \T^3} u_p(q)=u_{min}(p),p \in
\T^3
\end{equation}
and hence for all $ p\in U_\delta(0)$ the equality
\begin{equation}\label{q0}
 u_p(q_0(p))=\min_{q\in \T^3} u_p(q)=u_{min}(p)=u_{min}(-p)=u_{-p}(q_0(-p))=
 u_{p}(-q_0(-p))
\end{equation} holds.
For each $p \in U_\delta(0)$ the point $q_0(p)$ is the unique
non-degenerate minimum of the function $u_{p}(\cdot)$ and hence the
equality \eqref{q0} yields  $ q_0(-p)=-q_0(p),\,p\in U_\delta(0).$

The asymptotics \eqref{q0 asym} follows from the fact that
$q_0(\cdot)$ is an odd function and its coefficient
$-\frac{l_2}{l_1}$ is calculated using the identity
${\bigtriangledown } u(p, q_0(p))\equiv 0,\,p\in U_\delta(0).$

$(ii)$ The equality $u_{min}(p)=u_{p}(q_0(p))$ and  asymptotics
\eqref{q0 asym} yields asymptotics \eqref{min.raz}.

$(iii)$ Since the function $u_{\min}(\cdot)$ is even (see
\eqref{evenn}) we conclude that  if $q_0(p)\in \T^3$ is the minimum
point of $u_p(\cdot)$ then $-q_0(p)\in \T^3$ is the minimum point of
$u_{-p}(\cdot)$. Hence
 $q_0(-p)=-q_0(p).$
\end{proof}

Set
$$
\C_+=\{z\in \C: Re z>0\},\quad  \R_+=\{x\in \R: x>0\},\quad
\R_+^0=\R_+\cup \{0\}.
$$

Let $u(\cdot,\cdot)$ be the function defined on $U_{\delta}(0)\times
\T^3$ as
 \begin{equation*}\label{w}
\tilde u(p,q)=u_p(q+q_0(p))-u_{\min}(p).
\end{equation*}

For any $p\in {\T}^3$ we define an analytic function $D(p,\cdot)$ in
$\C_+$ by
\begin{equation*} D(p,w)= \int\limits_{{\T}^3}
\frac{\varphi^2(q+q_0(p))dq}{\tilde u(p,q)+w^2}.
\end{equation*}

\begin{lemma}\label{razl}  Assume Assumption \ref{hypoth} and \ref{hyp.varphi}.
 Then for any  $w\in \R_+^0$ the function
  $D(\cdot,w)$ is of class $C^{(2)}(U_{\delta}(0)),$
  the right-hand derivative of $D(0,\cdot)$ at $w=0$ exists  and
  the following decomposition
\begin{equation*}
D(p,w)=D(0,0)-\frac{4\sqrt{2}\pi^2 \varphi^2(0)}
  {l_1^{\frac{3}{2}}(det U)^{\frac{1}{2}}}
w+D^{res}(w)+D^{res}(p,w)
\end{equation*}
holds, where
 $D^{res}(w)=O(w^{1+\theta})$ as $ w\to +0$ and
 $D^{res}(p,w)=O(p^2)$ as $p\to 0$ uniformly in $w\in
\R_+^0.$
\end{lemma}

\begin{proof}  $(i)$ For any $p\in U_{\delta}(0)$ the
point $q_0(p)$  is the non-degenerate minimum of the function
$u_p(\cdot)$ (see Lemma \ref{minimum}) and $q_0(\cdot)\in
C^{(2)}(U_{\delta}(0))$. Since $u_{min}(\cdot)\in
C^{(2)}(U_{\delta}(0))$ by definition of $D(\cdot,\cdot)$ and
Assumptions \ref{hypoth} and \ref{hyp.varphi} we obtain that the
function $D(\cdot,w)$ is of class $C^{(2)}(U_{\delta}(0))$ for any
$w\in \R^0_+,$ where $C^{(n)}(U_{\delta}(0))$ can be defined in the
same way as $C^{(n)}(\T^3).$

One can easily see that
 \begin{align*}
u_p(q+q_0(p))-u_{\min}(p) =\frac{l_1}{2}(Uq,q)+o(|p||q|^2)+o(|q|^2)
\,\,as\,\,|p|,|q| \to 0.
\end{align*}
Therefore for some $C>0$ and all $w\in \R^0_+$ and $i,j=1,2,3$ the
inequalities hold
\begin{align}\label{Estimate1}
&\Big |\frac{\partial^2}{\partial p_i \partial p_j}
 \frac{\varphi^2(q+q_0(p))}{\tilde u(p,q)-w^2}\Big |\leq C|q|^{-2},\quad p,q\in U_\delta(0)
\end{align}
and
\begin{align}\label{Estimate2}
&\Big |\frac{\partial^2}{\partial p_i \partial p_j}
 \frac{\varphi^2(q+q_0(p))}{\tilde u(p,q)-w^2}\Big |\leq C,\quad p\in U_\delta(0), q\in \T^3
 \setminus U_\delta(0).
\end{align}

Lebesgue's dominated convergence theorem implies that
$$
\frac{\partial^{2}}{\partial p_i\partial p_j}D(p,0)=\lim_{w\to
0+}\frac{\partial^{2}}{\partial p_i\partial p_j}D(p,w),p\in
U_{\delta}(0).
$$

Repeatedly applying  Hadamard's lemma (see \cite{Zor} V.1, p. 512)
we obtain
\begin{align*}
D(p,w)=D(0,w)+\sum_{i=1}^3\frac{\partial}{\partial
p_i}D(0,w)p_i+\sum_{i,j=1}^{3}H_{ij}(p,w) p_i p_j,
\end{align*}
where for any $w\in \R^0_+$ the functions
$H_{ij}(\cdot,w),\,i,j=1,2,3$ are continuous in $U_\delta(0)$ and
\begin{equation*}
    H_{ij}(p,w)=\frac{1}{2}\int_{0}^1\int_{0}^1
    \frac{\partial^2}{\partial
p_i \partial p_j}D(x_1x_2p,w)dx_1dx_2.
\end{equation*}

Estimates \eqref{Estimate1} and \eqref{Estimate2} give
$$
|H_{i,j}(p,w)|\leq \frac{1}{2}\int_{0}^1\int_{0}^1
    |\frac{\partial^2}{\partial
p_i \partial p_j}D(x_1x_2p,w)|dx_1dx_2\leq C\big
(1+\int\limits_{U_\delta(0)}q^{-2}\varphi^2(q+q(p))dq \big )
$$
for any $p\in U_\delta(0)$ uniformly in $w\in \R_+^0.$

Since for any $w\in \R_+^0$ the function $D(\cdot,w)$ is even in
$U_\delta(0)$ we have
$$
\frac{\partial}{\partial p_i} D(p,w)\Big |_{p=0}=0,\quad i=1,2,3.
$$

$ii)$ Now we show that there exists a right-hand derivative of
$D(0,\cdot)$ at $w=0$ and for some $C>0$ the following relations
hold
\begin{equation}\label{partial D}
\lim_{w\to 0+}
  w^{-1}(D(0,w)-D(0,0))=
 \frac{2\sqrt{2}\pi^2\,\mu \varphi^2(0)}{l_1^{\frac{3}{2}}(det U)^{\frac{1}{2}}}  ,
\end{equation}
\begin{equation}\label{partial<}
\big |  \frac{\partial }{\partial w}D(0,w)- \frac{\partial
}{\partial w}D(0,0) \big |<C w^{\theta} ,\quad w\in \R^0_+.
\end{equation}

Indeed, the function $D(0,\cdot)$ can be represented as
\begin{equation*}
D( 0, w) =D_1(w)+D_2(w)
\end{equation*}
with
\begin{equation*}\label{I1}
D_1(w)=\int\limits_{U_{\delta}(0)}
\frac{\varphi^2(q)}{\tilde{u}(0,q)+w^2 }dq,\,w\in \C_+
\end{equation*}
and
$$
D_2(w)=\int\limits_{{\T}^3\setminus U_{\delta}(0)}
\frac{\varphi^2(q)}{\tilde{u}(0,q)+w^2 }dq,\,w\in \C_+.
$$

Since the function $\tilde{u}(0,\cdot)$ is continuous on the compact
set $\T^3\setminus U_\delta(0)$ and has a unique minimum at $q=0$
there exists $M_1>0$ such that $|\tilde{u}(0,q)|>M_1$ for all $q\in
\T^3\setminus U_{\delta}(0).$

Then by  $\varphi(\cdot)\in C^{(2)}(\T^3)$  we have
\begin{equation*}\label{D2}
|D_2(w)-D_2(0)|\leq C w^2,\,w\in \R_+^0
\end{equation*}
for some $C=C(\delta)>0.$

Now, let us consider
\begin{align}\label{D1}
D_1(w)-D_1(0)=\int\limits_{U_\delta(0)}
w^2[(\tilde{u}(0,q)+w^2)\tilde{u}(0,q)]^{-1}\varphi^2(q)dq .
\end{align}

The function  $\tilde{u}(0,\cdot)$ has a unique non-degenerate
minimum at $q=0.$ Therefore,  by virtue of the Morse lemma (see
\cite{Fed}) there exists a one-to-one mapping $q=\psi(t)$ of a
certain ball $u_\gamma(0)$  of radius $\gamma>0$ with the center at
\, $t=0$ to a neighborhood $\tilde W(0)$ of the point \, $q=0$ such
that
\begin{align}\label{eps=t2}
\tilde{u}(0,\psi(t))=t^2
\end{align}
with $\psi(0)=0$ and for the Jacobian $J_\psi(t) \in
\cB(\theta,U_\delta(0))$ of the mapping $q=\psi(t)$ the equality
holds
$$
J_\psi(0)=\sqrt{2} l_1^{-\frac{3}{2}}(det U)^{-\frac{1}{2}}.
$$

In the integral in \eqref{D1} making a change of variable
$q=\psi(t)$ and using the equality \eqref{eps=t2} we obtain

\begin{align}\label{D1-D1}
D_1(w)-D_1(0)=-\frac{w^2}{2}\int\limits_{u_\gamma(0)}
\frac{\varphi^2(\psi(t))J_\psi(t)}{t^2(t^2+w^2)}dt.
\end{align}

Going over in the integral in \eqref{D1-D1} to spherical coordinates
$t=r\omega,$ we reduce it to the form
\begin{equation*}\label{I2}
D_1(w)-D_1(0)=-\frac{w^2}{2} \int_0^{\gamma}\frac{F(r)}{r^2+w^2} dr,
\end{equation*}
with
$$
F(r)=\int_{\S^2}\varphi^2(\psi(r\omega))J_\psi(r\omega)d\omega,
$$
where $\S^2$ is the unit sphere in $\R^3$ and $d \omega$ is the
element of the unit sphere in this space.

 Using $\varphi\in C^{(2)}(\T^3)$ and $J_\psi\in
\cB(\theta,U_\delta(0))$ we see that
\begin{equation}\label{GEld}
|F(r)-  F(0)|\leq C r^{\theta},\,\, r\in [0,\delta].
\end{equation}

Applying the inequality \eqref{GEld} it is easy to see that
\begin{equation*}\label{F(0)}
\lim\limits_{w\to 0+}\frac {D_1(w)-D_1(0)}{w}=2\sqrt{2}\pi
l_1^{-\frac{3}{2}} \varphi^2(0)
  (det U)^{-\frac{1}{2}}.
\end{equation*}
Hence we have that there exists a right-hand derivative of
$D_1(\cdot)$ at $w=0$ and the equality \eqref{partial D} and the
inequality \eqref{partial<} hold.
\end{proof}
\section{Proof of the main results}
{\it\bf{The proof of Theorem \ref{pos.eig}}}. $(i)$ Part (i) resp.
part (ii) of Assumption \ref{Lambda} yields
$$
\Delta_{\mu_0}(p,u_{\min}(p))<\Delta_{\mu_0}(0,m)=0,\,0\neq p\in
\T^3.
$$
 resp.
$$\Delta_{\mu_0}(p,m)>\Delta_{\mu_0}(0,m)=0,\,0\neq p\in
\T^3.
$$
Since
 $ \lim_{z\to -\infty} \Delta_{\mu_0}(p,z)=1$ and
$\Delta_{\mu_0}(p,\cdot)$ is monotonously decreasing on
$(-\infty,u_{\min}(p))$ we conclude that the function
$\Delta_{\mu_0}(p,z)$ has an unique solution $e_{\mu_0}(p)$ in
$(m,u_{\min}(p)).$
 Lemma \ref{delta=0} completes the proof of Part $(i)$
of Theorem \ref{pos.eig}.

$(ii)$  Let $\mu>\mu_0.$ We have
\begin{equation}\label{d<D}
    \Delta_{\mu}(p,z)<\Delta_{\mu_0}(p,z)
\end{equation}
for all $ p\in \T^3,\,z\leq u_{\min}(p).$

Set $e_{\mu_0}(0)=m$. By assertion $(i)$ of Theorem \ref{pos.eig}
and Lemma \ref{delta=0} yield   the value of the function
$e_{\mu_0}(\cdot)$ at the point $p\in \T^3$ satisfies
 $m<e_{\mu_0}(p)<u_{\min}(p),$  $ p\neq 0,$  and   $m=e_{\mu_0}(0),p=0,$   and
$\Delta_{\mu_0}(p,e_{\mu_0}(p))=0.$

  By \eqref{d<D} we have
$\Delta_{\mu}(p,e_{\mu_0}(p))< \Delta_{\mu_0}(p,e_{\mu_0}(p))= 0,
p\in \T^3$ and hence by  Lemma \ref{delta=0} for any $p\in \T^3$
there exists the number
 $e_{\mu}(p)$ such that
\begin{equation*}
e_{\mu}(p)\in (-\infty,e_{\mu_0}(p))\quad \mbox{
 and}\quad \Delta_{\mu}(0,e_{\mu}(p))= 0.
\end{equation*}

 Hence  Lemma \ref{delta=0}   completes  the proof of Theorem
\ref{pos.eig}.\qed

{\it\bf{The proof of Theorem \ref{main}}} follows from Lemma
\ref{razl} if we take into account the equality
$\Delta_{\mu}(p,z)=1-\mu\Lambda(p,z)=1-\mu
D(p,\sqrt{u_{\min}(p)-z})$ and that $w=(u_{min}(p)-z)^{1/2}\geq 0$
for $z\leq u_{min}(p).$\qed

{\it\bf{The proof of Corollary \ref{razl.lemma.natijasi.}}} follows
from Theorem \ref{main} and Lemmas \ref{resonance}, \ref{zeroeigen}.
\qed

{\it\bf{Proof of Corollary \ref{D.ineq}}}.
 Let the
operator $h_{\mu_0}(0)$
 have a threshold  energy resonance then $\varphi(0)\neq 0$ (see Lemma
 \ref{resonance}).
 One has the asymptotics (see part
$(ii)$ of Lemma \ref{minimum})
\begin{equation*}\label{malfa}
u_{\min}(p)=m+(l_1l_2-l^2)(2l)^{-1}(Up,p)+o(|p|^2)\quad \mbox{as}
\quad p\to 0
\end{equation*}
and Corollary \ref{razl.lemma.natijasi.} yields \eqref{c<(.,.)<c}
for some positive numbers $c_1,c_2$.

 The positivity and continuity
of the function $\Delta_{\mu_0}(\cdot,m)$ on the compact set
$\T^3\setminus U_\delta(0)$ implies \eqref{(.,.)>c}. \qed

{\it\bf{Proof of Corollary \ref{fin.ineq}}}.
 By Lemma
\ref{zeroeigen} we have  $\varphi(0)=0$ and $\Delta_{\mu_0}(0,m)=0.$
Taking into account that  $\mu_0=\Lambda^{-1}(0,m)$, where the
function $\Lambda(\cdot,\cdot)$ is defined by \eqref{Lamb}, we get
the equality
$$
\Delta_{\mu_0}(p,m)=\mu_0(\Lambda(0,m)-\Lambda(p,m)).
$$
Then Assumption \ref{Lambda} completes the proof. \qed

\section{Appendix}
Here we show that there is some important  subclasses of the family
of Friedrichs models (see, e.g., \cite{ALMM,CL}) and for these
subclasses the assumptions in Section 2 are fulfilled.

Let \begin{align}\label{u form} \hat{u}(p,q)=\varepsilon (p)+
\varepsilon (p-q)+ \varepsilon (q),
\end{align} where $\varepsilon(p)$
is a real-valued
  conditionally negative definite function on $\T^3$ and hence
  \begin{itemize}
  \item[(i)] $\varepsilon $ is an even function,
  \item[(ii)] $\varepsilon(p)$ has a minimum at $p=0$.
  \end{itemize}

  Recall  that a complex-valued bounded
function $\varepsilon:\T^3\rightarrow \C$ is called conditionally
negative definite if $\varepsilon(p)=\overline{\varepsilon(-p)}$ and
\begin{equation}\label{nn}
  \sum_{i,j=1}^{n}\varepsilon(p_i-p_j)z_i\bar z_j\le 0
 \end{equation}
for any $n\in \N$, for all  $p_1, p_2, .., p_n\in \T^3$ and all
${\bf z}=(z_1, z_2, ..., z_n)\in \C^n$ satisfying
$\sum_{i=1}^nz_i=0$ (see, e.g., \cite{ALMM,RSIV}).

\begin{assumption}\label{end}
  Assume that $\varepsilon(\cdot)$ is a real-valued conditionally
negative definite  function on ${\T}^3$ having a unique
non-degenerate minimum at the origin and all third order partial
derivatives of $\varepsilon(\cdot)$ are continuous and belong to
$\cB (\theta,U_{\delta}(0)).$
\end{assumption}

\begin{lemma}\label{maximum}
Assume Assumption \ref{end} and $ u(p,q)=\hat{u}(p,q+q_0(p)).$ Then
part (i) of Assumption \ref{Lambda} is fulfilled.
\end{lemma}

\begin{proof}
$(i)$ By the definition of $\Lambda(p,u_{\min}(p))$ we have
\begin{align*}\label{L.Even}
&\Lambda(p,u_{\min}(p))-\Lambda(0,u_{\min}(0))\no\\
&=\int\limits_{{\T}^3}
\frac{2(u_0(t)-u_{\min}(0))-[u_p(t)+u_{-p}(t)-2u_{\min}(p)]}
{(u_p(t)-u_{\min}(p))(u_{-p}(t)-u_{\min}(p))(u_0(t)-u_{\min}(0)}\varphi^2(t)dt.
\end{align*}

According to $u_{p}(0)=u_{\min}(p)$ for any $0\neq p\in \T^3$ we
arrive to the inequality
$$u_0(q)-u_{min}(0)>\frac{u_p(q)+u_p(-q)}{2}-u_{ min}(p),\quad
\mbox{ a.e.,}\quad q\in \T^3,
$$
which is equivalent to the inequality in Lemma 5 \cite{ALMM} and
proves part $(i)$ of Assumption \ref{Lambda}.
\end{proof}
\begin{lemma}\label{maximum}
Let $  u(p,q)$ be of the form \eqref{u form}. Then part (ii) of
Assumption \ref{Lambda} is fulfilled.
\end{lemma}
\begin{proof}
 The real-valued (even) conditionally negative definite
function $\varepsilon$ admits the (L\'evy-Khinchin) representation
(see, e.g., \cite{BCR,ALMM})
\begin{equation}\label{L-Kh}
\varepsilon(p)=\varepsilon(0)+\sum_{s\in
\Z^3\setminus\{0\}}(\\cos(p,s)-1)\hat \varepsilon (s),\quad p\in
\T^3,
\end{equation}
which is equivalent to the requirement that the Fourier coefficients
$\hat \varepsilon(s)$ with $s\ne 0$   are non-positive, that is,
\begin{equation}\label{non-neg}
 \hat \varepsilon(s)\le 0, \quad s\ne 0,
\end{equation}
and the series $\sum_{s\in \Z^3\setminus\{0\}}\hat \varepsilon (s)$
converges absolutely.

Since $u$ and $|\varphi|$ are even the function $\Lambda(\cdot)$ is
also even. Hence the equality
\begin{equation*}
u_0(t)-\frac{u_p(t)+u_p(-t)}{2}= \sum_{s\in
\Z^3\setminus\{0\}}\hat{\varepsilon}(s)(1+\\cos (t,s))(1-\\cos
(p,s))
\end{equation*}
yields the representation
\begin{equation}\label{Lamb-lamb}
\Lambda(0,m)-\Lambda(p,m)= \frac{1}{2}\sum_{s\in
\Z^3\setminus\{0\}}(-\hat{\varepsilon}(s))(1-\\cos
(p,s))\int\limits_{{\T}^3} (1+\\cos(t,s))F(p,t)dt+\tilde{B}(p),
\end{equation}
where
$$F(p,\cdot)=\frac{[u_p(\cdot)+u_{-p}(\cdot)-2m]}
{(u_p(\cdot)-m)(u_{-p}(\cdot)-m)(u_0(\cdot)-m)}\varphi^2(\cdot)
$$
and
$$
\tilde{B}(p)=\frac{1}{4}\int\limits_{{\T}^3} \frac{
[u_p(t)-u_{-p}(t)]^2}
{(u_p(t)-m)(u_{-p}(t)-m)(u_0(t)-m)}\varphi^2(t)dt.
$$
Set
\begin{equation*}
B(p,s)=\int\limits_{{\T}^3} (1+\cos(t,s))F(p,t)dt.
\end{equation*}

 We rewrite the
function $B(p,s)$ as a sum of two functions
\begin{equation*}
B^{(1)}_{\delta}(p,s)= \int\limits_{{\T}^3\setminus U_\delta(0)}
(1+\cos(t,s))F(p,t)dt
\end{equation*}
and
\begin{equation*}
B^{(2)}_{\delta}(p,s)= \int\limits_{U_\delta(0)}
(1+\cos(t,s))F(p,t)dt.
\end{equation*}

Let $\chi_\delta(\cdot)$ be the characteristic function of $
U_\delta(0).$ Choose $\delta>0$ such that $$ mes\{ (\T^3\setminus
U_\delta(0))\cap supp\,\, \varphi\}>0 .$$

Set $F_{\delta}(p,\cdot)=(1-\chi_\delta(\cdot))F(p,\cdot).$ Then for
all $p\in \T^3$ and a.e. $$t\in (\T^3\setminus U_\delta(0))\cap
supp\,\,\varphi(\cdot)$$ the function $F_{\delta}(\cdot,\cdot)$ is
strictly positive.  Since the function $u$ has a unique minimum at
$(0,0)$ and
$\varphi$ is continuous on $\T^3$ we have that $
F_{\delta}(p,\cdot), p\in \T^3$ belongs to the Banach space
$L^1(\T^3).$ Then for some (sufficiently large) $R>0$ and
(sufficiently small) $c_1(\delta)>0$ and for all $|s|\leq R$, $p\in
\T^3$ we have the inequality
\begin{equation*}
B^{(1)}_{\delta}(p,s)= \int\limits_{{\T}^3}
(1+\cos(t,s))F_\delta(p,t)dt>c_1(\delta)>0.
\end{equation*}
 The
Riemann-Lebesgue lemma  yields
$$
B^{(1)}_{\delta}(p,s) \to
\int\limits_{{\T}^3}F_{\delta}(p,t)dt>0,\,p\in \T^3 \quad
\mbox{as}\quad s\to \infty.
$$
The continuity of the function \,\,$$\tilde F(p)=
\int\limits_{{\T}^3} F_{\delta}(p,t)dt $$ on the compact set $\T^3$
yields that for all $p\in \T^3$ and $|s|>R$ the inequality
$B^{(1)}_{\delta}(p,s)\geq c_2(\delta)$ holds.

Thus for $c(\delta)=\min\{c_1(\delta),c_2(\delta)\}$ the inequality
$B^{(1)}_{\delta}(p,s)\geq c$ holds for all $s\in \Z^3,\,p\in \T^3.$
So $B^{(2)}_{\delta}(p,s)\geq 0,\,s\in \Z^3,\,p\in \T^3$ yields
$B(p,s)>c,\,s\in \Z^3,\,p\in \T^3.$ Taking into account the
inequalities $\tilde{B}(p)\geq 0,\, p\in \T^3$ and
$\hat{\varepsilon}(s)\leq 0,\,s\in \Z^3\setminus\{0\}$ (see
\eqref{non-neg}) from the representations \eqref{L-Kh} and
\eqref{Lamb-lamb} we have
\begin{equation*}
\Lambda(0,m)-\Lambda(p,m)\geq c (\varepsilon(p)-\varepsilon(0)).
\end{equation*}
This together with the assumptions on $\varepsilon(\cdot)$ completes
the proof of Lemma \ref{maximum}.
\end{proof}
 {\bf Acknowledgement} The authors
would like to thank Prof.~R.~A.~Minlos and ~K.~A.~Makarov for
helpful discussions about
 the results of the paper. This work was
supported by the DFG 436 USB 113/6 Project and the Fundamental
Science Foundation of Uzbekistan. S.~N.~Lakaev and Z.~I.~Muminov
gratefully acknowledge the hospitality of the Institute of Applied
Mathematics and of the IZKS of the University of Bonn.We are greatly
indebted to the anonymous referee for useful comments.


\begin{thebibliography}{99}

\bibitem{AGH} Albeverio, S., Gesztesy F., and H{\o}egh-Krohn R.:
 The low energy expansion in non-relativistic scattering
theory. Ann. Inst. H. Poincar\'e Sect. A (N.S.) {\bf  37},
   1--28 (1982).

\bibitem{AGHH}
Albeverio S., Gesztesy F., H{\o}egh-Krohn R., and Holden H.:
 {\it Solvable Models in Quantum Mechanics}.
 Springer-Verlag, New York, 1988; 2nd ed. (with an appendix by P.~Exner),
 Chelsea, AMS, 2004.

\bibitem{ALMM}  {Albeverio ~S.~ ,  Lakaev ~S.~N.~,Makarov ~K.~A.~,
Muminov  ~Z.~I.~: { The Threshold Effects for the Two-particle
Hamiltonians on Lattices,}
 Comm.Math.Phys. {\bf 262}(2006), 91--115 .

\bibitem{ALzMahp}  { Albeverio ~S.~, Lakaev ~S.~N.~and
 Muminov ~Z.~I.~}: { Schr\"{o}dinger operators on lattices. The Efimov effect and
discrete spectrum asymptotics.} Ann. Henri Poincar\'{e}. {\bf 5},
(2004),743--772.


\bibitem{ALzMarX06}  { Albeverio ~S.~, Lakaev ~S.~N.~and
 Muminov ~Z.~I.~}:  The threshold effects for a family
 of Friedrichs  models under rank one
 perturbations.  arXiv:math.SP/0604277 v1 12 Apr 2006

\bibitem{ALZMarX06}  { Albeverio ~S.~, Lakaev ~S.~N.~and
 Muminov ~Z.~I.~}: 
On the number of eigenvalues of a model operator associated to a
system of three-particles on lattices. arXiv:math-ph/0508029 v2 14 Aug 2006

\bibitem{CL}
Carmona ~R. and Lacroix ~J.:
 {\it Spectral theory of random
Schr\"odinger operators. Probability and its Applications}, 1990,
Birkh\"auser Boston.

\bibitem{BCR}
Berg C.~, Christensen ~J.~P.~R.~, and  Ressel ~P.~, {\it Harmonic
analysis on semigroups. Theory of positive definite and related
functions.}
 Graduate Texts in Mathematics,
Springer-Verlag, New York, 1984. 289 pp.

\bibitem{J-K} Jensen ~A. and Kato ~T.:
 Spectral properties  of Schr\"odinger operators and
time-decay of the wave functions.
 Duke Math. J. {\bf 46}.
 583--611
(1979).

\bibitem{Ef}  {~Efimov ~V.~ }: Energy levels of three resonantly interacting
particles, Nucl. Phys. A {\bf 210} (1973), 157--158.



\bibitem{Fad}  { Faddeev ~ L.~ D.~}: On a model of Friedrichs in the theory
of perturbations of the continuous spectrum(Russian). Trudy Mat.
Inst. Steklov {\bf 73}(1964), 292--313.
\bibitem{Fed}  {M. ~V. ~Fedoryuk}: Asymptotics of integrals and
 series [in Russian], Nauka, Moscow (1987).

\bibitem{FIC} Faria da Veiga ~P.~A., Ioriatti ~L., and O'Carroll ~M.:
Energy-momentum spectrum of some two-particle lattice Schr\"odinger
Hamiltonians.
 Phys. Rev. E (3) {\bf  66}, \, 016130, 9 pp. (2002).
\bibitem{Fri} { Friedrichs ~K.~ O.}: On the perturbation of continuous spectra.
  Communications on Appl. Math.  {\bf 1}(1948), 361--406.

\bibitem{G-Sh} Graf ~G.~M. and Schenker ~D.:
$2$-magnon scattering in the Heisenberg model.
 Ann. Inst. H. Poincar\'e Phys. Th\'eor.
 {\bf 67},
  91--107 (1997).

\bibitem{K-S} Klaus M. and Simon B.:
   Coupling constants thresholds in
non-relativistic quantum mechanics. I. Short range two body case.
Ann. Phys. {\bf 130},
 251--281 (1980).

\bibitem{K-M}  Kondratiev Yu.~G. and R.~A.~Minlos R. A.:
 One-particle subspaces in the stochastic $XY$ model.
  J. Statist. Phys. {\bf 87},
   613--642 (1997).





\bibitem{Ltsp86}  {  Lakaev~S.~N.~}: Some spectral properties of the generalized
Friedrichs model, (Russian) Trudy Sem. Petrovsk. No. 11 (1986),
210--238, 246, 248; translation in J. Soviet Math. {\bf 45} (1989),
no. 6, 1540--1565.

\bibitem{Ltmf92}  { Lakaev~S.~N. ~}: {Bound states and resonances for the N-particle
discrete Schr\"{o}dinger operator,} Theor. Math. Phys. {\bf 91}
(1992), No.1, 362-372.

\bibitem{Lfa93}  {Lakaev~S.~N. ~}: {The Efimov's Effect of a system of Three
Identical Quantum lattice Particles,} Funkcionalnii analiz i ego
priloj. , {\bf 27} (1993), No.3, pp.15--28, translation in Funct.
Anal.Appl.

\bibitem{LMfa04}  {Lakshtanov~ E.~ L.~  and Minlos ~R.~ A.~ }: The spectrum of two-particle bound
states of transfer matrices of Gibbs fields (an isolated bound
state). (Russian) Funktsional. Anal. i Prilozhen. {\bf 38}(2004),
No.3, 52--69; (translation in Funct. Anal. Appl. 38 (2004), No. 3,
202--216)

\bibitem{LMfa05}  {Lakshtanov ~E. ~L.~   and ~Minlos ~R.~ A.~ }: The spectrum of two-particle bound
states of transfer matrices of Gibbs fields (fields on a
two-dimensional lattice: adjacent levels). (Russian) Funktsional.
Anal. i Prilozhen. {\bf 39} (2005), No. 1, 39--55, ( translation in
Funct. Anal. Appl. 39 (2005), no. 1, 31--45)

\bibitem{MalMin}  { Malishev ~V.~ A.~ and  Minlos ~R.~ A.~}: Linear infinite-particle
operators. Translations of Mathematical Monographs, 143. American
Mathematical Society, Providence, RI, 1995.



\bibitem{Mat} Mattis ~D.~C.:
The few-body problem on a lattice.
 Rev. Modern Phys. {\bf 58},
  361--379 (1986).

\bibitem{M-S} Minlos ~R.~A. and Suhov ~Y.~M.:
On the spectrum of the generator of an infinite system of
interacting diffusions.
 Comm. Math. Phys. {\bf 206},
  463--489 (1999).



\bibitem{Mog} Mogilner ~A.:
 Hamiltonians in solid state physics as
 multi-particle discrete Schr\"{o}dinger operators: Problems and
 results.
 Advances in Soviet Mathematics {\bf 5},
  139--194 (1991).



\bibitem{Rauch} Rauch ~J.:
 Perturbation theory for eigenvalues and resonances
of Schr\"{o}dinger Hamiltonians.
 J. Funct. Anal. {\bf 35},
  304--315 (1980).

\bibitem{RSIII}  Reed ~M. and  Simon ~B.:
{\it Methods of modern mathematical physics. III: Scattering
theory,}
 Academic Press, New York, 1979.


\bibitem{RSIV}  Reed ~M. and  Simon ~B.:
{\it Methods of modern mathematical physics. IV: Analysis of
Operators,}
 Academic Press, New York, 1979.



\bibitem{Sob}   Sobolev ~A.~V.:
The Efimov effect. Discrete spectrum asymptotics.
 Commun. Math. Phys. {\bf 156},
  127--168 (1993).

\bibitem{Tam94}  {Tamura ~H.~ }: Asymptotics for the number of negative eigenvalues of three-body
Schr\"odinger operators with Efimov effect. Spectral and scattering
theory and applications, 311--322, Adv. Stud. Pure Math., 23, Math.
Soc. Japan, Tokyo, 1994.

\bibitem{Wang}  { Wang ~X.~P.~}:  On the existence of the $N$- body Efimov effect, J.
Funct. Anal. {\bf 95} (1991), 433--459.




\bibitem{Yaf1}   Yafaev ~D.~R.:
 On the theory of the discrete spectrum of
the three-particle Schr\"{o}dinger operator.
 Math. USSR-Sb. {\bf 23},
  535--559  (1974).

\bibitem{Yaf2} Yafaev ~D.~R.:
 The virtual level of the Schr\"odinger equation.
J. Soviet. Math., {\bf 11},
 501--510 (1979).

\bibitem{Yaf3}  Yafaev ~D.~R.:
 {\it  Scattering theory: Some old and new problems},
 Lecture Notes in Mathematics, 1735.
 Springer-Verlag, Berlin, 2000, 169 pp.


\bibitem{Zh}  Zhizhina ~E.~A.:
 Two-particle spectrum of the generator for stochastic model
of planar rotators at
 high temperatures.
 J. Statist. Phys.  {\bf 91},
  343--368 (1998).


\bibitem{Zor}  {Zorich ~V.~ A.~ }:  Mathematical analysis I.
Springer-Verlag Berlin Heildelberg 2004.}


\end{thebibliography}
\end{document}